\input amstex
\documentstyle{amsppt}
\magnification=1200
\pagewidth{6.5 true in}
\pageheight{9.0 true in}
\nologo
\tolerance=10000
\overfullrule=0pt
\define\m{\bold m}
\define\ra{\rightarrow}
\define\R{R_{\geq Nq}}

\define\Jq{J^{[q]}}

\define\inc{\subseteq}

\define\Iq{I^{[q]}}
\define\x{\underline{x}}
\define\hgt{\text{height}}
\define\sat{^{\text{sat}}}

\topmatter
\title
The saturation of Frobenius powers of ideals 
\endtitle
\rightheadtext{Frobenius powers of ideals}
\author
Craig Huneke
\endauthor
\address
Department of Mathematics, University of Kansas,
Lawrence, KS 66045
\endaddress
\email
huneke\@math.ukans.edu
\endemail
\thanks
The author 
thanks the NSF for
partial support. I thank Yongwei Yao for many valuable comments.
\endthanks
\abstract
Some special cases are proved in which localization commutes with
tight closure. 
\endabstract
\endtopmatter
 
\document
\centerline{Dedicated to Robin Hartshorne on the occasion of his 60th birthday}
\bigskip
\head 1. Introduction
\endhead

The work of Robin Hartshorne includes many papers which analyzed and charted
the structure of local cohomology modules. Among these papers were those
which studied local cohomology in positive characteristic. 
In 1977, Robin Hartshorne and Robert Speiser \cite{HS} wrote a ground-breaking paper
on the action of Frobenius on local cohomology. Now over twenty years later,
this Frobenius action, while still not completely understood, plays a pivotal and
crucial role in the theory of tight closure, especially in its
relationship to the classification of singularities via the action of
Frobenius (cf. \cite{Sm2, Ha, HaW, MS, Wa1,2}). One of the frustrating problems in the theory of tight closure
has been our inability to decide whether tight closure commutes with
localization. This question is closely related to understanding
the local cohomology of Frobenius powers of modules. In fact localization would
be proved if one could show that a certain set of associated primes is finite,
and prove either condition $(LC^*)$ or $(LC)$ as in Definition 1.2 below.
See \cite{Ka1, Hu1, Chapter 12}. Condition $(LC)$ concerns the saturation
of the Frobenius powers of an ideal. 

We begin by stating relevant definitions. For more references and information about
tight closure we refer the reader to \cite{HH1, Hu1}, and for  survey articles
on this topic to \cite{Br,Sm1}. 
The phrase ``characteristic $p$'' always means positive and prime characteristic $p$.
We will use `$q$' throughout this paper to denote a variable power of the
characteristic $p$.
 Tight closure
is a method which requires  reduction to characteristic $p$, although
the theory has now been developed for any Noetherian ring containing a field.
 
We set $R^o$ to be the set of all elements of $R$ not in any minimal prime
of $R$.  The definition of tight closure for ideals is:

\definition{Definition 1.1} Let $R$ be a Noetherian ring of characteristic $p >
0$. Let $I$
be an ideal of $R$. An element $x\in R$ is said to be in the tight closure of $I$ if
there exists an element $c\in R^o$ such that for all
large $q= p^e$, $cx^q\in I^{[q]}$, where $I^{[q]}$ is the ideal generated by the $q$th
powers of all elements of $I$.
\enddefinition

One of the main open problems about tight closure is whether it
commutes with localization. New impetus has been given to study this
problem in the graded case by the work of Gennady Lyubeznik and Karen
Smith \cite{LS}.  Studying the localization problem leads naturally to
the following definition:

\definition{Definition 1.2} Let $R$ be a graded ring over a field $R_0$ of
positive characteristic $p$ with irrelevant ideal $\m = \oplus_{i\geq 1} R_i$.
Set $R_{\geq N}$ to be the ideal generated by all forms of degree at least $N$.
 Let $I$ be a homogeneous ideal of $R$. We say that
\it $I$ satisfies condition $(LC)$ \rm if there is a constant $N$ such that for all
$q$, $$\R\cdot H^0_{\m}(R/\Iq) = 0.$$ We say that \it $I$ satisfies condition $(LC^*)$ \rm
if there is a constant $N$ such that for all
$q$, $$\R\cdot H^0_{\m}(R/(\Iq)^*) = 0.$$
\enddefinition

We can similarly define these conditions for a local Noetherian ring
$(R,\m)$, by replacing $\R$ by $\m^{Nq}$ in the definitions. 

Observe that we can rephrase these conditions in terms of the saturation
of an ideal. Let $R$ be as in Definition 1.2, and let $I$ be an ideal of
$R$. The \it saturation \rm of $I$ (with respect to $\m$), denoted $I\sat$,
is the ideal $I\sat = \{r\in R|\,\, \m^nr\inc I\, \text{for some}\, $n$\}$. With
this definition, $I$ satisfies $(LC)$ iff there exists a constant $N$ 
such that for all
$q$, $\R\cdot(\Iq)\sat\inc \Iq$, and satisfies condition $(LC^*)$ iff
there exists a constant $N$  
such that for all 
$q$, $\R\cdot((\Iq)^*)\sat\inc (\Iq)^*$.

A priori there is no relationship between the conditions $(LC)$ and $(LC^*)$.
If $R$ has an isolated singularity then the theory of test elements \cite{HH2}
proves there is an $\m$-primary ideal $\tau$ such that $\tau\cdot(\Iq)^*\inc \Iq$
for all $q$ and all ideals $I$, and in this case it is obvious that 
$I$ satisfies $(LC)$ iff it satisfies $(LC)^*$. In general, however, we do
not know the precise relationship between the two properties.
As far as we know all homogeneous ideals $I$ satisfy both of these two conditions.
These conditions play a critical role in the question of whether the tight closure
of $I$ commutes with localization \cite{AHH, Kat1,2}. The paper \cite{Hu2} will
give necessary and sufficient conditions for tight closure to localize which
involve slightly weaker conditions.

In this paper we will prove a condition relating whether a homogeneous
ideal $I = (J,y)$ satisfies condition $(LC)$ (respectively $(LC^*)$) to
whether $J$ satisfies condition $(LC)$ (respectively $(LC^*)$). This is
done when $J:y$ has large enough codimension. We are then able to use
this theorem to give several new cases in which tight closure commutes
with localization. See Corollary 3.2, Corollary 3.7, and Corollary 3.8.
Corollary 3.2 recovers a recent result of Adela Vraciu \cite{V} with a much
different proof.

\remark{Remark 1.3} One fact we will use without further comment is that if $R$ is graded and
$I$ is homogeneous, then the tight closure of $I$ is also homogeneous,
\cite{HH3, (4.2)}. In particular if $I$ is homogeneous, then the 
associated primes of $(\Iq)^*$ are homogeneous for all $q$.

We will also need to use the existence of test elements. From \cite{HH3, (3.6) and (4.2)}
it follows that if $R$ is graded there exists a homogeneous element $c\in R^o$ such
that for all ideals $I$ and all $x\in I^*$, $cx^q\in \Iq$ for all $q >> 0$. If $R$ is further assumed to
be a domain, then we can choose $c\ne 0$ and homogeneous such that $x\in I^*$
iff for all $q$, $cx^q\in \Iq$.
\endremark
\bigskip
\head 2. The Main Theorem 
\endhead
\bigskip
In this section we prove our main result, which is a criterion for an ideal
$I = J + Ry$ to satisfy either $(LC)$ or $(LC^*)$ in terms
of properties of $J$ and $J:y$. We will use this criterion in the next section to give
our main applications. Before stating the theorem, some elementary lemmas
on graded rings are helpful.

\proclaim{Lemma 2.1} Let $(R,\m)$ be a Noetherian graded ring over a field $R_0 = K$.
There exists a constant $L$ such that for all nonnegative integers $N,M$,
$$R_{\geq(N+M+L)}\inc R_{\geq N}\cdot R_{\geq M}.$$
In particular, for all $n \geq L,\,  R_{\geq(N+M+1)n}\inc R_{\geq Nn}\cdot R_{\geq Mn}$.
In addition for all $n$,  $R_{\geq(N+M+L)n}\inc R_{\geq Nn}\cdot R_{\geq Mn}$
\endproclaim

\demo{Proof}. By choosing a system of parameters of the same degree $d$,
one sees that $R$ is module finite over the subring $A$ generated by these
parameters. Observe that $A$ is generated by forms of the same degree $d$ so
that $A_{sd} = (A_d)^s$. 
Say the module generators of $R$ over $A$ can be chosen to be of degrees at most $k$. We claim that
$L = k + 2d$ satisfies the statement of the Lemma. Consider the graded
piece $R_{N+M+2d+k+i}$ for some $i\geq 0$. Set $m = N+M+2d+k+i$. We can write
this piece as $\oplus_s R_{m-sd}A_{sd}$ for integers $s$ such that $m-sd\leq k$.
Set $p = \lceil \frac N d \rceil$ and $q = \lceil \frac M d \rceil$
Then $N+M\leq pd + qd < N+M+2d$. Hence $sd\geq m-k\geq N+M+2d > pd + qd$,
so that $s\geq p+q$. Then $R_{m-sd}A_{sd} = R_{m-sd}(A_d)^s\inc
 R_{\geq (m-sd)}(A_d)^p(A_d)^q\inc R_{\geq N}\cdot R_{\geq M}$. 
The last two statements follow immediately.\qed
\enddemo

\proclaim{Lemma 2.2} Let $(R,\m)$ be a Noetherian graded ring over a field $R_0 = K$ of
positive characteristic $p$. Let $J$ be an arbitrary homogeneous ideal and $z$ a form of
degree $k$ such that $\sqrt{(J,z)} = \m$. Then there exist   constants
$N_1,N_2$ such that for all $q = p^e$ and all integers $l$,
 $$R_{\geq(N_1q + N_2lkq)}\inc J^{[q]} + (z^{N_2lq}).$$
\endproclaim

\demo{Proof} 
Choose homogeneous generators $x_1,...,x_n$ of $J$, and choose powers $x_1^{a_1},...,
x_n^{a_n}$ and $z^{N_2}$ such that the degree of all of these elements is
the same, namely $N_2k$. Clearly it suffices 
replace the $x_i$ by these powers and $z$ by $z^{N_2}$, call their
common degree $d$,  and then prove
there exists a single constant $N$ such that 
$$ R_{\geq({Nq + ldq})}\inc J^{[q]} + (z^{lq}).$$ 

Let $A = K[x_1,...,x_n,z]$ be the graded subring of $R$ generated by
the given elements. 
 Because $(x_1,...,x_n,z)$ is $\m$-primary,
$R$ is a finite $A$-module and there exist a finite set of generators up
to degree $L$. Consider $R_{M+L}$. This piece is contained in 
$\sum_{i=0}^L R_iA_{M+L-i}\inc A_MR$. Set $M = nqd +qld$. Then we
see that $R_{\geq(nqd+lqd+L)}\inc A_{\geq(nqd+lqd)}R = (x_1,...,x_n,z)^{nq+lq}R\inc
(x_1,...,x_n)^{nq}R + z^{lq}R\inc J^{[q]} +  z^{lq}R$. Setting
$N_1 = nd+L$ then gives the statement of the Lemma.\qed\enddemo

\bigskip
\proclaim{Theorem 2.3} Let $(R,\m)$ be a Noetherian graded ring of dimension $d$
 over a field $R_0 = K$ of
positive characteristic $p$. Let $J = (x_1,...,x_n)$ be a homogeneous ideal
and let $y$ be a form of degree $e$. Set  $I = J + Ry$. Assume that $\hgt(J:y)\geq d-1$.
\roster
\item If $J$ satisfies $(LC)$ then $I$ satisfies $(LC)$. 
\item If $J$ satisfies $(LC^*)$, then $I$ satisfies $(LC^*)$.
\endroster
\endproclaim

\demo{Proof} We begin with (1).  Choose $z\in \m$ general of degree $k$. 
For some fixed $N$, $R_{\geq Nq}\cdot (\Jq)\sat\inc \Jq$ for all $q$. 
For a fixed $q$ and all large $L$, $(\Iq)\sat = \Iq:z^{Lq}$. Consider 
the short exact sequence,
$$0\ra ((\Jq:z^{Lq}) + Ry^q)/\Iq\ra (\Iq:z^{Lq})/\Iq\ra (\Iq:z^{Lq})/((\Jq:z^{Lq}) + Ry^q)\ra
 0\tag{2.4}$$

The first term in this short exact sequence is killed by $R_{\geq Nq}$ by our
assumption, since $(\Jq:z^{Lq})\inc (\Jq)\sat$. Consider the last term. It is
isomorphic with 
$$
((\Jq, z^{Lq}):y^q)/((\Jq:y^q)+Rz^{Lq}).\tag{2.5}
$$
This last isomorphism follows by moding out $\Jq$ and applying the fact that
in general for elements $a,b$ in a ring there is an isomorphism 
$$(a:b)/((0:b)+(a))\cong (b:a)/((0:a)+(b)).$$

Since $(J:y)^{[q]}\inc (\Jq:y^q)$, we know that $(J:y)^{[q]} + Rz^{Lq}$ is contained
in the denominator ideal of (2.5). By Lemma 2.2 there exists
constants $N_1, N_2$ such that for all $l$ and all $q = p^e$,
$$
R_{\geq(N_1q+N_2lkq)}\inc (J:y)^{[q]} + Rz^{N_2lq}.$$
We claim that for all $l$ and all $q$, $R_{\geq(N_1 +e)q}$ annihilates the module (2.5),
with $L = N_2l$.
Suppose that $r\in (\Jq, z^{N_2lq}):y^q$ is homogeneous. If the degree of $r$
is at least $(N_2lk-e)q$, then $R_{\geq(N_1 +e)q}r\inc R_{\geq(N_1q+N_2lkq)}\inc
 (J:y)^{[q]} + Rz^{N_2lq}\inc ((\Jq:y^q)+Rz^{N_2lq})$
and so the coset of $r$ in (2.5) is $0$. If deg$(r) < (N_2l-e)q$, then writing
$ry^q = w + vz^{N_2lq}$ with $w\in \Jq$  we see that $v$ may be chosen to be $0$ since the
degree of the term $vz^{N_2lq}$ is larger than the degree of $ry^q$. Then
$ry^q\in \Jq$ and again the coset of $r$ is $0$. 

It follows that $R_{\geq(N_1+e)q}\cdot R_{\geq Nq}$ annihilates
 the middle module
of the sequence (2.4) for all $l$, with $L = N_2l$. By choosing $l >> 0$, 
we can identify the middle module with $(\Iq)\sat/\Iq$. By Lemma 2.1,
there exists another constant $K$ such that
 $R_{\geq Kq}\inc R_{\geq(N_1+k)q}\cdot R_{\geq Nq}$ for all $q$, which
 finishes the proof of the (1).

The proof of (2) is similar but somewhat tricky. 
Let $N$ be chosen so that $R_{\geq Nq}((\Jq)^*)\sat\inc (\Jq)^*$ for all
$q$, and let $N'$ be chosen so that for all $q$ $R_{\geq (N'q + lqk)}\inc (J:y)^{[q]} + Rz^{lq}$
for some suitable power $z$ (of degree $k$) of a general form. This is possible by
Lemma 2.2 and the assumption that the height of $J:y$ is  at least $d-1$.  

We claim that $R_{\geq Nq}\cdot R_{\geq (N'q+eq)}((\Iq)^*)\sat\inc (\Iq)^*$.
To prove this claim, let $r\in ((\Iq)^*)\sat, u\in  R_{\geq (N'q+eq)}$,
and $v\in R_{\geq Nq}$ be arbitrary homogeneous elements. Let $z$ be as in the above paragraph.
For some $l >> 0$, $z^{lq}r\in (I^{[q]})^*$. By Remark 1.3 there exists a 
homogeneous element $c\in R^o$ such
that for all $q' >> 0$ there exist  homogeneous elements $s_{qq'}$ with
$$cr^{q'}z^{lqq'}\in J^{[qq']} + s_{qq'}y^{qq'}.$$

First suppose that the term $s_{qq'}y^{qq'}$ does not appear for infinitely
many $q'$. Then $cr^{q'}z^{lqq'}\in J^{[qq']}$ for all such $q'$ which implies that
$rz^{lq}\in (\Jq)^*$ \cite{Hu1, Exercise 2.2} (this is the local case, but
the proof given in the Appendix of \cite{Hu1} applies equally well to the graded case).  Then $r\in ((\Jq)^*)\sat$, and $vr\in (\Jq)^*\inc (\Iq)^*$.
A fortiori $uvr\in (\Iq)^*$ which proves our claim. We may therefore assume that
the term $s_{qq'}y^{qq'}$ appears for all large $q'$.
In this case 
deg$(s_{qq'}) = \text{deg}(cr^{q'}z^{lqq'}) - \text{deg}(y^{qq'}) =
\text{deg}(c) + q'\text{deg}(r) + lqq'k - eqq'$. If we multiply by
$u^{q'}$ we obtain that deg$(u^{q'}s_{qq'})\geq N'qq' + lqq'k$, and hence
$u^{q'}s_{qq'}\in (J:y)^{[qq']} + Rz^{lqq'}$. It follows that $u^{q'}s_{qq'}y^{qq'}\in
J^{[qq']} + Rz^{lqq'}y^{qq'}$, and hence that
 $$cr^{q'}z^{lqq'}u^{q'}\in J^{[qq']} + Rz^{lqq'}y^{qq'}.$$

Factoring out the common $z^{lqq'}$ we obtain that $cr^{q'}u^{q'} - t_{qq'}y^{qq'}\in 
(J^{[qq']})\sat$ for some element $t_{qq'}$.  By choice of $N$ it then holds that
$v^{q'}(cr^{q'}u^{q'} - t_{qq'}y^{qq'})\in (J^{[qq']})^*$ and hence there exists an
element $d\in R^o$ such that
$dv^{q'}(cr^{q'}u^{q'} - t_{qq'}y^{qq'})\in J^{[qq']}$. This holds for all large $q'$.
From this equation, $cd(ruv)^{q'}\in J^{[qq']} + Ry^{qq'}$ which forces
$ruv\in (\Iq)^*$. This proves our claim.

By Lemma 2.2 there exists a constant $L$ such that
$R_{\geq (N+N'+e+L)q}((\Iq)^*)\sat\inc (\Iq)^*$ for all
$q$ which proves (2) and finishes the proof of the Theorem.  \qed
\enddemo
\bigskip
\head 3. Application to the Localization of Tight Closure
\endhead
\bigskip

In order to apply the main theorem of the above section, one needs to find
an ideal $J$ satisfying either $(LC)$ or $(LC^*)$ to begin an induction. One
such class is given in our first proposition.

\proclaim{Proposition 3.1} Let $(R,\m)$ be a Noetherian graded equidimensional ring of
dimension $n$ over a field $R_0 = K$ of
positive characteristic $p$.  Let $J = (x_1,...,x_{n-1})$ be a homogeneous ideal
generated by homogeneous parameters $x_i$ (so that the height of $J$ is $n-1$).
Then $J$ satisfies both $(LC)$ and $(LC^*)$.
\endproclaim
 
\demo{Proof} The fact that $J$ satisfies $(LC^*)$ comes from the fact that
$(J^{[q]})^*$ is unmixed (see \cite{AHH,(5.20)}).
For the second, choose a form $z$ general. Set $x_n = z$.
Consider the the Koszul homology $H_1(x_1,...,x_n;R)$. Since this homology
is killed by the $\m$-primary ideal generated by the $x_i$, it has finite length.
This module is graded. If $(r_1,...,r_n)$ represents a homogeneous element of
degree $N$ in $H_1(x_1,...,x_n;R)$ (so that $\sum r_ix_i = 0$), then
deg$(r_i) + d_i = N$, where deg$(x_i) = d_i$. Since the module has finite
length, 
there is an integer $N$ such that $[H_1(x_1,...,x_n;R)]_i = 0$ for all
$i\geq N$. Using \cite{AH, Proposition 4.10} we obtain that for all $q = p^e$ and
for all $l$,
$$[H_1(x_1^q,...,x_{n-1}^q,x_n^{lq};R)]_i = 0$$
for all $i\geq N + qD + lqd_n$, where $D = \sum_{i=0}^{n-1} d_i$.
 
Let $r\in (\Jq)^{sat}$ be homogeneous of degree $M$.
 For some $l >> 0$, $x_n^{lq}r\in \Jq$; the power $l$
a priori depends upon $q$. If we write $x_n^{lq}r = \sum_{i=0}^{n-1} r_ix_i^q$,
with each term homogeneous of the same degree,
then the element $(r_1,...,r_{n-1},-r)$ represents an element of
 $H_1(x_1^q,...,,x_{n-1}^q,x_n^{lq};R)$
of degree $M + lqd_n$.  Suppose that $M\geq N+qD$. Then this element is trivial
in $H_1(x_1^q,...,,x_{n-1}^q,x_n^{lq};R)$ and this implies that $r\in J^{[q]}$.
 It follows that
$R_{\geq (D+N)q}\cdot r\inc \Jq$, which finishes the proof of the Proposition. \qed
\enddemo

Proposition 3.1 and Theorem 2.3 have several interesting corollaries. The first is due to
Adela Vraciu \cite{V} using very different methods:

\proclaim{Corollary 3.2} Let $(R,\m)$ be a Noetherian graded ring over a field $R_0 = K$ of
positive characteristic $p$. Let $I$ be a homogeneous ideal with dim$(R/I) = 1$.
\roster
\item If $R$ is equidimensional,
then $I$ satisfies both $(LC)^*$ and $(LC)$.
\item $(I^*)_W = (I_W)^*$ for
every multiplicatively closed set $W$ of $R$.
\endroster
\endproclaim

\demo{Proof} We use induction on $n = \mu(I)$. If $n = \hgt(I)$, then $I$ is generated
by part of a system of parameters. In this case  Proposition 3.1 gives that
$I$ satisfies both conditions. If $n > \hgt(I)$ we can choose minimal generators
$y_1,...,y_{n-1}, y_n = y$ for $I$ such that the height of $J = (y_1,...,y_{n-1})$
is still $\hgt(I)$. In this case $J$ and $y$ satisfy the conditions of Theorem 2.3
and since by induction $J$ satisfies both $(LC)^*$ and $(LC)$, so does $I$.

To prove the last statement, we may reduce to the case in which $R$ is a domain by \cite{AHH, (3.8)}.
Assume that the last statement does not hold. If after moding out by a minimal prime
of $R$ the ideal $I$ becomes $\m$-primary, we are done by \cite{HH1, (4.14)}. Henceforth
we assume this is not the case.
From \cite{AHH, (3.5)} it follows that
for some prime ideal $P$ containing $I$, $(I_P)^* \ne (I^*)_P$.
First assume that  $P$ is minimal over $I$. Choose an
element $t\notin P$ such that $t$ is in every other minimal primary 
component of $I$.
Suppose that $x\in (I_P)^*$.
By Remark 1.3, one can find a homogeneous element $c\in R^o$ and elements $w_q\notin P$ such
that $cw_qx^q\in \Iq$ for all $q$. 
The ideal $\Iq:cx^qt^q$ is then
$\m$-primary since it is homogeneous and not contained in any minimal prime
of $\Iq$.
Hence $cx^qt^q\in (\Iq)\sat$ for
all $q$. 
Since $I$ satisfies condition $(LC)$ 
there is an $N$ with $\R(cx^qt^q)\inc \Iq$ and we may assume by increasing $N$
that $R_N\nsubseteq P$.
Choosing $u\in R_N, u\notin P$ then gives $c(uxt)^q\in \Iq$ for all $q$ which
forces $utx\in I^*$. But then $x\in (I^*)_P$ as required.

Next suppose that P is a maximal ideal of $R$ which is not the irrelevant ideal.
Let $P_1,...,P_k$ be the minimal primes of $I$ contained in $P$ and let
$Q_1,...,Q_l$ be the minimal primes of $I$ which are not contained in $P$.
Choose $t\notin P$ homogeneous such that $t\in I_{Q_i}$ for all $i = 1,...,l$. This is
possible by the definition of the $Q_i$. Let $z\in (I_P)^*$. Using (1.3) there exist
elements $w_q\notin P$ and a homogeneous $c\ne 0$ such that $cw_qz^q\in \Iq$ for all
$q$. Consider the ideal $(\Iq:c(tz)^q)$. This ideal contains $w_q$ and also
contains a power of $P_1\cap ...\cap P_k$. It follows that this ideal is
not contained in any of the minimal primes over $I$ and so its radical is
an intersection of maximal ideals $\m, M_1,...,M_l$. There exists an integer
$N(q)$ such that $(\m M_1M_2...M_l)^{N(q)}c(tz)^q\inc \Iq$. But $\Iq$ is homogeneous
and all its associated primes are homogeneous. Thus $M_1M_2...M_l$ is not
contained in an arbitrary associated prime of $\Iq$ which forces
$\m^{N(q)}c(tz)^q\inc \Iq$.  
Hence $c(tz)^q\in (\Iq)\sat$ and there exists a constant
$N$ such that $\R(c(tz)^q)\in \Iq$ and we may choose $N$ large enough so that
$R_N\nsubseteq P$. Choose an element $u\in R_N$ such that
$u\notin P$. Then $u^q\in (\Iq:c(tz)^q)$ so that $c(utz)^q\in \Iq$. It follows
that $z\in (I^*)_P$ as $ut\notin P$.

Finally suppose that $P = \m$. Since every element not in $\m$ is a nonzerodivisor
on $\Iq$, it easily follows that $(I^*)_{\m} = (I_{\m})^*$ (see \cite{AHH, (3.5c)}. \qed\enddemo

It is well-known to experts that (3.2) gives the implication `$R$ is weakly
F-regular' $\implies$ `$R$ is F-regular'. We sketch a proof here. Recall
weakly F-regular means all ideals are tightly closed, while F-regular means
all ideals are tightly closed in every localization. Choosing an ideal $I$
maximal in $R$ and homogeneous such that $I_W$ is not tightly closed in
some localization, one can easily prove that $I$ must be an irreducible ideal
with the dimension of $R/I$ equal to one. This uses that all ideals are
tightly closed in $R$. Then apply the last statement of (3.2). 
Murthy \cite{Hu1, (12.2)} has shown that weakly F-regular rings which are
finitely generated over an uncountable field of positive characteristic are
F-regular.
Recently, Lyubeznik and Smith \cite{LS} have shown the stronger implication
that graded weakly F-regular rings are even strongly F-regular.

To extend  Corollary 3.2 to ideals of smaller height, we need a condition to
transfer the two conditions $(LC)^*$ and $(LC)$ from one ideal to one of lower
height. This is possible in some cases:

\proclaim{Proposition 3.3} Let $R$ be an equidimensional Noetherian graded ring
 over a field $R_0$  of characteristic $p > 0$ with irrelevant ideal $\m$.
Let $I$ be an ideal and assume that $I_P$ is generated by $g = \text{height}(I)$
elements for every minimal prime $P$ of $I$ of height $g$. Let $x_1,...,x_g$ be arbitrary set
of  $g$ elements in $I$ such that $\hgt(\x) = g$ and $(\x)_P = I_P$ for every minimal
prime $P$ of $I$ of height $g$.
\roster
\item  If $(LC^*)$ holds for
the ideal $J = (\x:I) + I$, then it holds for $I$.
\item If $(LC)$ holds for $J = (\x:I) + I$ then $(LC)$ also
holds for $I$ if in addition $R$ is Cohen-Macaulay.
\endroster
\endproclaim
 
Before beginning the proof, we need a crucial lemma.
 
\proclaim{Lemma 3.4} Let $R$, $I$ and $\x$ be as above.
\roster
\item $(\x^{[q]}:\Iq)\cap ((\Iq)^*)^{\text{sat}}\inc (\x^{[q]})^*$.
\item If $R$ is Cohen-Macaulay
then $(\x^{[q]}:\Iq)\cap (\Iq)\sat = (\x^{[q]})$.
\endroster
\endproclaim
 
\demo{Proof} We first prove (1).
Let $u\in (\x^{[q]}:\Iq)\cap ((\Iq)^*)^{\text{sat}}$ and choose $w\in \m$ a general form
of sufficiently high degree so that $wu\in (\x^{[q]}:\Iq)\cap (\Iq)^*$. Then
there exists an element $c\in R^o$ such that for all large $q'$,
$c(wu)^{q'}\in (\x^{[q]}:\Iq)^{[q']}\cap I^{[qq']}$. As
$(\x^{[q]}:\Iq)^{[q']}\inc (\x^{[qq']}:I^{[qq']})$ we obtain that
$c(wu)^{q'}((\x^{[qq']}:I^{[qq']}) + I^{[qq']})\inc \x^{[qq']}$. 
 By \cite{AHH, (5.20)} the ideal $(\x^{[qq']})^*$ is
unmixed; every associated prime has height $g$.
Since $(\x)_P = I_P$ for every minimal prime of $I$ of height $g$, the ideal
$((\x^{[qq']}:I^{[qq']}) + I^{[qq']})$ has height at least $g+1$. 
Then we must have that $c(wu)^{q'}\in (\x^{[qq']})^*$. It follows that
$wu\in  (\x^{[q]})^*$, and then the general choice of $w$ gives that
$u\in  (\x^{[q]})^*$ as required. 

The proof of (2) is similar but easier. If $R$ is Cohen-Macaulay, then
$(\x^{[q]})$ is unmixed and the desired equality follows immediately from
checking it after localization at the minimal primes above $\x$.
\enddemo
 
\demo{Proof of (3.3)} Let $I$ and $J$ be as in (3.3), and assume that
$r\in ((\Iq)^*)\sat$ in case (1) and in $(\Iq)\sat$ in case (2).
 Fix an integer $N$ such that  in case (1), $m^{Nq}((\Jq)^*)\sat\inc
(\Jq)^*$, and in case (2), $m^{Nq}((\Jq))\sat\inc \Jq$. 

In case (1) we know that
$$m^{Nq}r\inc (\Jq)^*\cap ((\Iq)^*)\sat,$$ while in case (2),
$$m^{Nq}r\inc \Jq\cap (\Iq)\sat.$$
Case (2) follows immediately from Lemma 3.4.2 since
$\Jq\cap (\Iq)\sat\inc (\x:I)^{[q]}\cap (\Iq)\sat + \Iq\inc
(\x^{[q]}:\Iq)\cap (\Iq)\sat + \Iq = \Iq$ by (3.4.2).
 
To prove case (1), 
choose $z$ an arbitrary element of $\m^N$. Then for all $q$,
$z^qr\in (\Jq)^*$, and so for all $q'$, $$cz^{qq'}r^{q'}\in (\x:I)^{[qq']} +I^{[qq']}.$$
As $r\in ((\Iq)^*)\sat$, there is a power
$N(q)$ of $\m$ such that $\m^{N(q)}r\inc (\Iq)^*$, which implies that
for all $q'$, $c (\m^{N(q)})^{[q']}r^{q'}\inc I^{[qq']}$, i.e.
$cr^{q'}\in (I^{[qq']})\sat$. Putting this together with the displayed
containment above, we see that
$$cz^{qq'}r^{q'}\in (\x:I)^{[qq']}\cap (I^{[qq']})\sat + I^{[qq']}.$$
Hence
$$cz^{qq'}r^{q'}\in (\x^{[qq']}:I^{[qq']})\cap (I^{[qq']})\sat + I^{[qq']}\inc (I^{[qq']})^*$$
by Lemma 3.4.2. This equation holds for all $q'$, proving
that $z^qr\in (\Iq)^*$. This means that $(\m^N)^{[q]}r\inc
(\Iq)^*$. If the number of generators of $\m^N$ is $l$, then it follows
that $\m^{Nlq}r\inc (\Iq)^*$. As $r$ was arbitrary in $((\Iq)^*)\sat$
we have proved that $\m^{Nlq}((\Iq)^*)\sat)\inc (\Iq)^*$, proving (1).\qed\enddemo

In a similar vein we can prove:

\proclaim{Proposition 3.5} 
Let $R$ be an equidimensional Noetherian graded ring
 over a field $R_0$  of characteristic $p > 0$ with irrelevant ideal $\m$.
Let $I$ be a homogeneous ideal and assume that $I_P$ is generated by $g = \text{height}(I)$
elements for every minimal prime $P$ of $I$ of height $g$. Let $x_1,...,x_g$ be arbitrary set
of  $g$ homogeneous elements in $I$ such that $\hgt(\x) = g$ and $(\x)_P = I_P$ for every minimal
prime $P$ of $I$ of height $g$.
Set $J = (\x:I) + I$. If $(J^*)_W = (J_W)^*$ for all multiplicatively closed
sets $W$ in $R$, then $(I^*)_W = (I_W)^*$ for all multiplicatively closed
sets $W$ in $R$.
\endproclaim

\demo{Proof} It is enough to prove this proposition for multiplicatively closed
sets of the form $W = R - Q$, where $Q$ is a prime ideal containing $I$ (see
\cite{AHH, (3.5)}). Let $z\in R$ be such that $z\in (I_Q)^*$. If $(\x:I)\nsubseteq Q$,
then $z\in (\x_Q)^* = (\x^*)_Q$ \cite{AHH,(8.1)} and so $z\in (I^*)_Q$ as required.
Hence we may assume that $J\inc Q$. Then $z\in (J_Q)^* = (J^*)_Q$. By multiplying
$z$ by an appropriate element not in $Q$ we may assume that $z\in J^*$. Let
$\{P_1,...,P_k\}$ be the minimal primes of $\x$ which are contained in $Q$,
and let $\{Q_1,...,Q_l\}$ be the minimal primes of $\x$ not contained in $Q$.
We can choose an element $t\in Q_1\cap...\cap Q_l$, $t\notin Q$ such that
$t$ is in $(\x)_{Q_i}$ for all $i$.  

Since $z\in J^*$, there exists an element $c\in R^o$ such that 
$cz^q\in (\x:I)^{[q]} + \Iq$ for all large $q$. Write $cz^q = a_q + b_q$
with $a_q\in (\x:I)^{[q]}$ and $b_q\in \Iq$. Since $z\in (I_Q)^*$, there
exists elements $w_q\notin Q$ and $d\in R^o$ such that $dw_qz^q\in \Iq$ for all $q >> 0$.
Then $dw_qa_q\in (\x:I)^{[q]}\cap \Iq\inc ((\x)^{[q]})^*$ by (3.4.1).
By \cite{AHH,(5.20)} $((\x)^{[q]})^*$ is unmixed. Since $w_q\notin Q$, it also holds
that $w_q\notin P_j$ for all $j$ and thus that $da_q\in (((\x)^{[q]})^*)_{P_j}$
for all $j$. The choice of $t$ then proves that $dt^qa_q\in (((\x)^{[q]})^*)_{P}$
for every minimal prime $P$ of $\x$. The unmixedness of $((\x)^{[q]})^*$
then proves that $dt^qa_q\in ((\x)^{[q]})^*$ for all $q >> 0$. Hence
$$cdt^qz^q = dt^qa_q + dt^qb_q\in ((\x)^{[q]})^* + \Iq\inc (\Iq)^*$$
for all $q >> 0$. It follows that $tz\in I^*$, and since $t\notin Q$, that
$z\in (I^*)_Q$ as desired. \qed\enddemo

\remark{Remark 3.6} Let $R$ be graded over an infinite field $R_0$ and suppose
that $I$ is a homogeneous ideal, $P_1,...,P_n$ are primes ideals and
$I_{P_j}$ is generated by $g$-elements for each $P_j$. Further assume that
for all $P_j\ne \m$, $R_N\nsubseteq P_j$ for all $N >> 0$ (e.g. if $R = R_0[R_1]$).
 Then we can find
$g$ homogeneous elements $x_1,...,x_g\in I$ such that $(x_1,...,x_g)_{P_j} = I_{P_j}$
for every $P_j$.

To see this claim, write $I = (y_1,...,y_k)$ with all $y_i$ homogeneous. If one of the
$P_j$ is the irrelevant ideal $\m$, there is nothing to prove as $I$ will be
generated by $g$ homogeneous elements. Suppose none of the $P_j$ is $\m$. 
For all large $N$ and all $1\leq j \leq n$, $R_N\nsubseteq P_j$, so we may choose
$z_j\notin P_1\cup P_2\cup ...\cup P_n$ homogeneous such that
each $z_jy_j$ is homogeneous of the same degree, and we still have
 $I_{P_j} = (z_1y_1,...,z_ky_k)_{P_j}$ for all $1\leq j \leq n$. We have
reduced to the case in which all $y_i$ have the same degree. But now any
$g$ general linear combinations of the generators with field elements will
generate $I$ locally at each $P_j$, and the resulting elements will be
homogeneous.

\endremark
\proclaim{Corollary 3.7}
Let $R$ be an equidimensional Noetherian graded ring
 over an infinite field $R_0$  of characteristic $p > 0$ with irrelevant ideal $\m$.
Let $I$ be an ideal with dim$(R/I) = 2$ and assume that $I_P$ is generated by $g = \text{height}(I)$
elements for every minimal prime $P$ of $I$ of height $g$.
Assume that $R_N\nsubseteq P$ for all minimal primes $P$ containing $I$ with
dim$(R/P) = 2$ and all $N >> 0$ (e.g. if $R = R_0[R_1]$).
Then $I$ satisfies $(LC^*)$. Furthermore, $(I^*)_W = (I_W)^*$ for every
multiplicatively closed set $W$ in $R$.
\endproclaim

\demo{Proof}
Let $x_1,...,x_g$ be arbitrary set
of  $g$ homogeneous elements in $I$ such that $\hgt(\x) = g$ and $(\x)_P = I_P$ for every minimal
prime $P$ of $I$ of height $g$. Such elements exist using Remark 3.6.
Set $J = (\x:I) + I$. Using (3.3.1) it suffices to prove that $J$ satisfies
$(LC^*)$. But dim$(R/J) < \text{dim}(R/I)$ and hence dim$(R/J)\leq 1$. 
Corollary 3.2 gives us that $(J^*)_W = (J_W)^*$ for every multiplicatively
closed set $W$ and $J$ satisfies $(LC^*)$. Then Proposition 3.5 proves that
 $(I^*)_W = (I_W)^*$
for every multiplicatively
closed set $W$ and $I$ satisfies $(LC^*)$. \qed\enddemo

We can generalize Corollary 3.7 to ideals of smaller height, but to do
so we need to assume the base ring is Cohen-Macaulay:

\proclaim{Corollary 3.8} Let $(R,\m)$ be 
a Cohen-Macaulay Noetherian graded ring
over an infinite field $R_0
= K$ of positive characteristic $p$. Assume that $R = R_0[R_1]$.  Let $I$ be a homogeneous ideal
such that $I_P$
is generated by $g = \hgt(I)$ elements for every homogeneous prime ideal $P$
such that dim$(R/P)\geq 2$. Then $I$ satisfies $(LC^*)$ and $(LC)$.
Furthermore, $(I^*)_W = (I_W)^*$ for every
multiplicatively closed set $W$ in $R$.

\endproclaim

\demo{Proof}
We use induction on dim$(R/I)$ to prove the first two statements. If dim$(R/I) \leq 1$, Corollary
3.2 gives the result. Henceforth we assume dim$(R/I)\geq 2$. Then
$I_P$ is generated by $g = \hgt(I)$ for every minimal prime $P$ of $I$ of height $g$.
Choose $\underline{x}$ as in Proposition 3.3, using Remark 3.6.
Using this proposition we see that it suffices to prove that
$J = (\underline{x}:I) + I$ satisfies $(LC^*)$ or $(LC)$ respectively.
 But our assumptions on $I$
prove that dim$(R/J) < \text{dim}(R/I)$. 
Then we are done using induction provided we prove that
$J_Q$ is generated by $g+1 = \hgt(J)$ elements for every prime $Q$ of $J$
such that dim$(R/Q)\geq 2$.
Let $Q$ contain $J$. If dim$(R/Q)\geq 2$, we know that $I_Q$ is generated
by $g$ elements, say $y_1,...,y_g$, which we may assume are homogeneous by
Remark 3.6. But in $R_Q$ we can then calculate
$J_Q$ as the image of $\underline{x}:\underline{y} + (\underline{y})$, and this
is simply $(\underline{y}) + R\Delta$ where $\Delta$ is the determinant of
the $g$ by $g$ matrix over $R_Q$ expressing the $x_i$ in terms of $y_i$.
Hence $J_Q$ is generated by $g+1$ elements. 
The induction finishes the proof.

To prove  $(I^*)_W = (I_W)^*$ for every 
multiplicatively closed set $W$ in $R$, we also use induction on dim$(R/I)$.
If dim$(R/I)\leq 1$, we can use Corollary 3.2. Define $J$ as in the above
paragraph. By the work above, dim$(R/J) < \text{dim}(R/I)$, and 
$J_Q$ is generated by $g+1 = \hgt(J)$ elements for every minimal prime $Q$ of $J$
such that dim$(R/Q)\geq 2$. By induction, $(J^*)_W = (J_W)^*$ for every 
multiplicatively closed set $W$ in $R$, and then by Proposition 3.5, we
can conclude $(I^*)_W = (I_W)^*$ for every 
multiplicatively closed set $W$ in $R$. \qed\enddemo

\bigskip
\centerline{\bf Bibliography}
\bigskip
\refstyle{A}

\enddocument